\documentclass{birkmult}

\usepackage{amsmath}
\usepackage{amssymb}
\usepackage{invol}

\begin{document}

%involtop.tex
\title[Involutions in Degenerate Type Groups]
{Involutions in groups of finite {M}orley rank of degenerate type}
%\vskip .1 true in November, 2006}	%DATE- omitted in final draft

\author[A.~Borovik]{Alexandre Borovik}
\address{School of Mathematics\\
University of Manchester\\
Oxford Road, Manchester M13 9PL\\U.K.}
\email{borovik@manchester.ac.uk} 

\author[J.~Burdges]{Jeffrey Burdges}
\address{Universit\'e de Lyon\\Universit\'e Lyon 1\\
CNRS, UMR 5208 Institut Camille Jordan\\
43, blvd du 11 novembre 1918\\
F-69622 Villeurbanne Cedex\\France}
\email{burdges@math.rutgers.edu}

\author[G.~Cherlin]{Gregory Cherlin}
\address{Department of Mathematics\\
Rutgers University\\
110 Frelinghuysen Rd.\\
Piscataway, NJ 08854\\U.S.A.}
%\email{OMIT} %For editorial use only: cherlin@math.rutgers.edu

\thanks{
\indent Second author supported by NSF Graduate Research Fellowship,
and DFG at Wurzburg and Bielefeld, grant Te 242/3-1,
and acknowledges the hospitality of the Universities of Birmingham and
Manchester.\\
\indent Third author supported by NSF Grant DMS-0100794.\\
\indent All authors thank the Newton Institute, Cambridge, for its
hospitality during the Model Theory and Algebra program, where the bulk
of this work was carried out, as well as CIRM for its hospitality at
the September 2004 meeting on {\em Groups, Geometry and Logic}, where
the seed was planted. Thanks to Alt{\i}nel for continued discussions
all along the way.
}

\begin{abstract}
In a connected group of finite Morley rank,
if the Sylow 2-subgroups are finite then they are trivial.
The proof involves a combination of model theoretic ideas with a
device originating in black box group theory.  
\end{abstract}


% \keywords{}

\maketitle
% \thispagestyle{empty}

%\markboth{{\sc Authors Here}}
%{{\sc Title here}}
%\pagestyle{myheadings}

%degenB01.tex

\section{Introduction}

Modern model theory can be viewed as a subject obsessed with notions
of dimension, with the key examples furnished by linear dimension on
the one hand, and the dimension of an algebraic variety (or, from
another point of view, transcendence degree) on the other.
There are several rigorous, and not always equivalent, notions of
abstract dimension in use. For historical reasons the one we
use is generally referred to as Morley rank. In the applications of
model theory, it is important that this dimension may be ordinal valued,
but the case of finite dimension continues to stand out. For example,
in the model theoretic approach to the Manin kernel in an abelian variety,
one enriches the underlying algebraically closed field with a
differential field structure, at which point the abelian variety
becomes infinite dimensional, but the Manin kernel itself is finite
dimensional, which accounts for a certain number of its fundamental
properties. 

For some time it was hoped that one would be able to classify
the ``one-dimensional'' objects arising in  model theory explicitly
and in complete generality, a hope which was dashed by a construction
of Hrushovski. But in diophantine applications, even after leaving the
algebraic category, one has in addition to the dimension notion a
topology reminiscent of the Zariski topology, and some 
very strong axioms, given in \cite{HZ-ZG}.
In this case one gets the desired algebraicity result, in
nondegenerate cases, with substantial diophantine
applications (cf.~\cite{Bo-ML,Sc-DG}).
\goodbreak

The work reported here concerns the Algebraicity Conjecture
(Cherlin/Zilber) which states that a simple group of finite Morley rank 
should be algebraic. This occupies a middle position between
the known results used in diophantine applications (where the focus in
any case is on abelian groups) and the more
ambitious conjectures which have been refuted. There is no assumption
of a topological nature, and the axioms are only those which occur in
general model theory, but considerably refined by results holding
in the specific context of groups, where the group action introduces a
degree of uniformity into the picture. 

This algebraicity conjecture also occupies a kind of middle position between
algebraic group theory and finite group theory. The identification of
the simple algebraic groups as Chevalley groups can be carried out
with relentless efficiency by examining maximal tori and their actions
on unipotent subgroups, thereby quickly revealing the associated root
system and giving the structure of the Weyl group. The classification
of simple finite groups is similar in outcome: setting aside the
alternating groups and a sackful of sporadic groups, one has some
sort of twisted Chevalley group, which can be identified by determining
the associated building, though the process is so intricate that by
the time the building is actually visible the whole group is equally visible.
The question arises whether either of these two approaches offers
anything for our more general problem. In particular cases, both do,
and they can even be combined. But this requires a certain supply of
elements of order two to implement.

Some time ago a project was launched to apply the
techniques of finite simple group theory in combination with relevant
notions of algebraic group theory, toward an analysis of hypothetical
nonalgebraic simple groups of finite Morley rank {\it containing
involutions,} which aims at pinning down the critical (minimal) 
configurations with as much precision as possible.  An early sketch of the
possibilities is found in \cite{BN}. 
In some cases the hypothesis of minimality is superfluous: for
example, if the group contains a nontrivial unipotent $2$-subgroup in
an appropriate sense \cite{ABC:book}.
Another such case is the case in point in the
present article.
We will prove the following.

\begin{theorem}\label{thm:main}
Let $G$ be a connected group of finite Morley rank
whose Sylow $2$-subgroup is finite. 
Then $G$ contains no involutions---the Sylow $2$-subgroup is trivial.
\end{theorem}

As a simple group will be either
finite or connected, this result tells us in the simple case that
the only exceptions are the finite simple groups.
The following version is a little sharper.

\begin{theorem}\label{thm:main2}
Let $G$ be a connected group of finite Morley rank, and $a\in G$ an
involution. Then the Sylow $2$-subgroup of $C(a)$ is
infinite.
\end{theorem}

This casts considerable light on the general program of determining
the possible $2$-Sylow structure in a hypothetical counterexample to
the Cherlin/Zilber Algebraicity Conjecture, according to which simple
groups of finite Morley rank should be Chevalley groups.

Groups of finite Morley rank are abstract groups equipped with 
a notion of dimension which assigns to every definable set $X$
a dimension, called {\em Morley rank} and denoted $\rk(X)$,
satisfying well known and fairly rudimentary axioms given for example 
in \cite{BN,Po-GrS}. Examples are furnished by algebraic groups over
algebraically closed fields, with $\rk(X)$ equal to the dimension
of the Zariski closure of $X$.
The Algebraicity Conjecture amounts to
the assertion that simple algebraic groups can be characterized, as
abstract groups, by the presence of a dimension. A striking feature of
this conjecture is the complete absence of any {\em topological}
assumptions; as far as this circle of ideas is concerned,
it is the main conjecture which survives without imposing topological
conditions. 

As the category of groups of finite Morley
rank is closed under finite direct products, it is easy to
fabricate non-algebraic examples, but to construct ``seriously''
nonalgebraic examples is a challenge, which so far has been best met
by Baudisch in \cite{Ba-UCG} with a nilpotent example.

To put the ongoing work in context, we need to present an
overall framework for the analysis, as well as the current status
of the  program.
One first defines {\it Sylow $2$-subgroups} in the natural way, and
one also considers the {\em connected components} of Sylow $2$-subgroups,
which are called Sylow\oo\ $2$-subgroups. The structure of these 
Sylow\oo\ $2$-subgroups is as
follows, in the context of groups of finite Morley rank.
$$\hbox{$S=U*T$\qquad $U$ $2$-unipotent, $T$ a $2$-torus}$$
That is, $S$ is the central product of groups $U$ and $T$ where 
$U$ is definable and connected, and of bounded exponent, while
$T$ is a divisible abelian $2$-group, not in general definable.
This simple fact provides a framework for further analysis.

It is known that necessarily $S=U$ or $S=T$ when the ambient group $G$
is simple, the other factor being trivial \cite{ABC:book}. 
Furthermore,
when $U$ is nontrivial the group $G$ is in fact algebraic.
These are substantial results relying on an analysis which would be
considered long by most standards, though unbelievably rapid by the
standards of finite simple group theory. So attention focusses on
the case $S=T$, which in fact is two cases:
$$\hbox{$S=T>1$ (odd type); \quad $S=T=1$ (degenerate type)}$$
We note that another formulation of the degeneracy condition is that
the full Sylow $2$-subgroup of the ambient group $G$ is {\em finite},
and that is the condition we have adopted in formulating our main
result above.

It is in odd type that the theory still relies on the assumption 
that the ambient simple group is a $K^*$-group; that is,
all its proper definable simple sections are algebraic. Given that, if
the ambient simple group is nonalgebraic then 
the group $T$ satisfies severe restrictions: it is known that its Pr\"ufer rank
(which is the same as the $2$-rank $m_2(T)$) is at most two,
which is a condition analogous to Lie rank at most two in the
algebraic case. A good deal more is known about this case, and a good
deal remains to be done here. 

Throughout, the case of degenerate type has occupied a peculiar
position. It has certainly not been ignored, but on the other hand
nothing approaching a systematic plan or point of view has
ever been found for dealing with this case. The elimination of
involutions from these groups has been arrived at unexpectedly,
and from an unusual direction: techniques used in the computational
theory of so-called {\em black box groups} provide a key ingredient in
the proof, and the remaining ingredients are purely model theoretic.
What we borrow is taken from one particular chapter in the theory 
which owes little to the
conventional analysis of finite simple groups, other than a
preoccupation with centralizers of involutions.
The sort of problem we are dealing with here (which is related to the
$Z^*$-theorem in finite group theory) would normally be approached
using character theoretic methods and transfer arguments, neither of
which have analogs in our setting. The black box methods become very
global and direct in our setting, having to do with the ranks
of fibers of appropriate covariant maps (first in a rudimentary
way in \S3, then more precisely in \S5).

Tuna Alt{\i}nel has pointed out  
that with Theorem \ref{thm:main} in hand, our subject
reaches a level of maturity sufficient to 
dispose of two ``chestnuts'' from the
general theory of groups of finite Morley rank by reduction to the
simple case and application of what amounts to a tripartite theory
(degenerate type, odd type, and {\em even type}---the last being the
case $S=U>1$). Here the fact that odd type groups are always studied
in an inductive setting could be an obstacle, but as it happens a
certain amount of information of a noninductive character has also
been obtained recently \cite{Ch-GT}.
The two results in question are the following.

\begin{proposition}\label{prop:centralizer}
Let $G$ be a connected and nontrivial group of finite Morley rank.
Then the centralizer of any element of $G$ is infinite.
\end{proposition}

\begin{proposition}\label{prop:evenorder}
Let $G$ be a connected group of finite Morley rank containing a
definable generic subset whose elements are of 
order $2^k$ for some fixed $k$.
Then $G$ has exponent $2^k$.
\end{proposition}

Both of these problems have been notoriously open.
The point of the second one, stressed by Poizat, is that
the corresponding result for algebraic groups is trivial,
 in view of the presence of a suitable topology \cite[p.~145]{Po-EG,Po-GrS}.

One can combine these two results as follows.

\begin{proposition}\label{prop:genericexponent}
Let $G$ be a connected group of finite Morley rank containing a
definable generic subset whose elements are of 
order $n$ for some fixed $n$. 
Then the Sylow $2$-subgroup $S$ of $G$ is unipotent, $G=S*C_G(S)$ is a
central product,
and $G/S$ is a group without involutions 
whose elements are generically of order
$n_0$, where $n_0$ is the odd part of $n$.
\end{proposition}

Thus for all practical purposes the study of the ``generic equation''
$x^n=1$ reduces to the case of $n$ odd. The analysis in that case is
complicated by the possible existence of simple $p$-groups of finite
Morley rank.

While Proposition \ref{prop:genericexponent} is expressed as a strong
form of Proposition \ref{prop:evenorder}, it also includes
Proposition \ref{prop:centralizer}: 
if $a\in G$ has a finite centralizer, then $a$ has finite order and 
its conjugacy class $a^G$ is generic in $G$.
Then applying Proposition \ref{prop:genericexponent}, evidently $G>S$
and thus $\bar G=G/S$ is a group without involutions in which the conjugacy
class of the image $\bar a$ of $a$ is still generic. But the same
applies
to $\bar a^{-1}$ and thus $\bar a$, $\bar a^{-1}$ are conjugate, and
this produces involutions in $G/S$, a contradiction.

We remark that these last results (and Theorems \ref{thm:p-odd},
\ref{thm:p-odd2} following)  
rely on a considerable body of
material, not all of it fully published at this time, and the
reader may prefer to consider them as conditional---specifically,
conditional on the classification of groups of even type
\cite{ABC:book}.

We will generalize Theorem \ref{thm:main} as follows.

\begin{theorem}\label{thm:p-odd}
Let $G$ be a connected group of finite Morley rank
whose Sylow $p$-subgroups are finite. Then $G$ contains no elements of
order $p$.
\end{theorem}

However the proof for odd $p$ involves a reduction to the statement
for $p=2$, though, as will be seen, parts of the two proofs can be done
uniformly. 

Then we can continue on to the parallel to Theorem \ref{thm:main2}.

\begin{theorem}\label{thm:p-odd2}
Let $G$ be a connected group of finite Morley rank, and $a\in G$ a
$p$-element. Then the Sylow $p$-subgroup of $C(a)$ is
infinite.
\end{theorem}

Returning to the case $p=2$, we will also show the following.

\begin{theorem}\label{thm:generation}
Let $V$ be a $4$-group acting definably on a connected group $H$ of
finite Morley rank and degenerate type.
Then 
$$H=\<C_H\oo(v):v\in V^\times\>$$
\end{theorem}

In fact one first applies Theorem \ref{thm:main} to reduce 
to the involution-free case (and in that setting, connectivity is no
longer needed). However this more general form is
the form one would actually want, and may actually be useful, eventually,
in emancipating the odd type analysis from the $K^*$-hypothesis,
a line of development which remains to be explored. In this line, the
following is completely open.

\begin{problem}
Can one show that a simple group of finite Morley rank and degenerate type
has no nontrivial involutory automorphism?
\end{problem}

Evidently the hypothesis of simplicity is necessary, but is not easy to
bring to bear on the problem.

Finally, we should mention an outstanding open question of a general
nature concerning groups of finite Morley rank.

\begin{conjecture}[Genericity Conjecture]\label{con:genericity}
Let $G$ be a connected group of finite Morley rank.
Then the union of the connected definable 
nilpotent subgroups of $G$ contains a definable generic subset of $G$.
\end{conjecture}

This conjecture reflects quite well a useful property of connected
algebraic groups. Its truth would greatly simplify 
the proof of our main result. Conversely, 
some of our methods might contribute to the general
analysis of this problem. Jaligot has shown \cite{Ja-Ge}
that statements closely related to the Genericity Conjecture 
have strong implications
for the general structure theory of groups of finite Morley rank, and
mesh well with the existing Carter subgroup theory.

\medskip

Notation follows \cite{BN}. We mention the notation $d(g)$ for the
{\em definable closure} of $g$ in the {\em group theoretic} sense: this is
the smallest definable subgroup containing $g$, and coincides with
$d(\< g\>)$. See \cite[\S5 {\it et passim}]{BN} for this.

We thank Alt\i nel for a number of useful discussions, not confined to
the points already mentioned above.

% End of file degenB01.tex
 %Intro
%degenB02.tex

\section{Approximating $d(g)$}\label{sec:technical}

The present section has a purely technical character.
We will make extensive use of the function $d(g)$ which assigns to an
element $g$ the smallest definable subgroup containing it,
which we think of as the definable group ``generated by'' $g$.
Unfortunately the function $d$ is not in general definable
(in an algebraic geometric context the parallel remark would be that
the family $\{d(g):g\in G\}$ is not an algebraic family). Indeed,
for $g$ of finite order the group $d(g)$ is the finite cyclic group
generated by $g$, which is typically of unbounded order, and hence
cannot be given in a uniformly definable way.

We will replace $d$ by a definable approximation $\hat d$ with
essentially the same properties from a practical point of view.
There is no ``canonical'' approximation, but we will record all the
properties of $\hat d$ used in the present paper, which we think
provides a robust approximation to the function $d$.

A word on the terminology used in the next lemma is in order.
We are interested in
$p$-torsion (elements of order $p^n$) for all primes $p$, with special
interest in the case $p=2$. We will work generally with the hypothesis
that Sylow $p$-subgroups have finite exponent. 
There are 
two distinct notions of Sylow $p$-subgroup, both reasonably natural, 
in this context: either a
maximal $p$-subgroup, as usual, or else a maximal locally nilpotent
$p$-subgroup (local nilpotence can be 
replaced here by anything similar---local finiteness or solvability for
example). 

With either definition, there is at present 
no {\it Sylow theory} for $p$ odd, except in the framework of solvable groups. 
Note however that the following conditions are equivalent.

\begin{enumerate}
\item
 Every abelian  $p$-subgroup of $G$ has finite exponent.
\item
 Every Sylow $p$-subgroup of $G$ has finite exponent.
\item
The $p$-torsion of $G$ has bounded exponent.
\end{enumerate}

Since these are successively stronger conditions, it suffices
to show that the first implies the last. If the $p$-torsion
of $G$ has unbounded exponent, then in a  saturated elementary
extension $G^*$ of $G$ there will be an infinite quasicyclic $p$-subgroup
$A$. Then $d(A)$ will be a definable abelian $p$-divisible group
containing an element of order $p$. The existence of such a subgroup
passes from $G^*$ to $G$ and contradicts $(1)$.

Since all of these conditions are equivalent, we may express them as
follows: ``$G$ has Sylow $p$-subgroups of bounded exponent.''
The reader may possibly prefer to keep clause $(3)$ in mind.

The practical effect of this condition is the following, which may be
surprising if one thinks in terms of ``cyclic'' subgroups;
the correct intuition is furnished by
{\em Zariski closures} of cyclic subgroups.

\begin{lemma}\label{L:d}
Let $G$ be a group of finite Morley rank, $p$ a prime. Suppose
that $G$ has Sylow $p$-subgroups of bounded exponent, and let $p^n$ be
a bound on the exponent. Then for any element
$g\in G$, the following hold.
\begin{itemize}
\item $d(g^{p^n})$ is uniquely $p$-divisible.
\item $d(g)$ has a cyclic Sylow $p$-subgroup $S$.
\item $d(g)=d(g^{p^n})\times S$.
\end{itemize}
\end{lemma}
\begin{proof}
Let $A=d(g)$, a definable abelian group. 
One has $A=A_0\oplus S$ with $A$ $p$-divisible and $S$ a $p$-group of
finite exponent. 
By our hypothesis, $A_0$ has no $p$-torsion and thus $S$
is the Sylow $p$-subgroup of $A$. We can write $g=as$ with $a\in A_0$,
$s\in S$ and then $g\in A_0\times \<s\>$, so $d(g)=A_0\times \<s\>$.
All that needs to be checked at this point is that $A_0=d(g^{p^n})$.

We have $g^{p^n}=a^{p^n}$. Since every definable subgroup of $A_0$ is
$p$-torsion free and thus uniquely $p$-divisible, it follows that
$a\in d(g^{p^n})$ and thus $g\in d(g^{p^n})\times S$. Hence
$d(g)=d(g^{p^n})\times S$ and $d(g^{p^n})=A_0$.
\end{proof}

The picture provided by the foregoing lemma should be borne in mind
throughout. 
Note that tori behave quite differently; their generic elements are
dense, and thus the Sylow $p$-subgroups of $d(g)$ in general can be
very far from cyclic. 

\begin{lemma}
Let $p$ be a prime, and 
let $G$ be a group of finite Morley rank with Sylow $p$-subgroups of
bounded exponent. Then there is a definable function $\hat d(a)$,
from elements of $G$ to definable subgroups of $G$, with the following
properties.
\begin{enumerate}
\item $d(a)\le \hat d(a)$;
\item If $d(a)=d(b)$, then $\hat d(a)=\hat d(b)$;
\item For $g\in G$, we have $\hat d(a^g)=\hat d(a)^g$;
\item $\hat d(a)$ is abelian;
\item The groups $d(a)$ and $\hat d(a)$ have the same Sylow
  $p$-subgroup;
\item If $x\in G$ conjugates $a$ to its inverse, then $x$ normalizes
  $\hat d(a)$ and acts on it by inversion.
\end{enumerate}
\end{lemma}
\begin{proof}
Consider the following two functions.
\begin{itemize}
\item $d_1(a)=Z(C(a))$. 
\item $d_2(a)=d_1(a)^q\<a\>$ where $q$ is a bound on the
order of the $p$-torsion in $G$.
\end{itemize}

One sees easily that $d_1$ is definable and 
satisfies our first four conditions (cf.~\cite[\S5.1]{BN}).
As $[d_2(a):d_1(a)^q]\le q$ it follows easily that $d_2$ is also a
definable function, and it also satisfies condition $(1)$ and inherits
conditions $(2-4)$ from $d_1$. 

Furthermore, $d_2(a)$ also satisfies the fifth condition,
since $d_1(a)$ is abelian and $d_1(a)^q$ is $q$-torsion free.

Now to achieve the final point, let $d_3(a)$ be the subgroup of 
$d_2(a)$ consisting of elements inverted by every element that inverts
$a$. 
\end{proof}

The following is in a similar vein.

\begin{lemma}
Let $G$ be a group of finite Morley rank, and $p$ a prime, with
Sylow $p$-subgroups of $G$ of bounded exponent.
Let  $H$ be a definable subgroup of $G$,
and $a\in N(H)$ of order $p$ modulo $H$.
Then $d(a)\intersect a H$ contains a $p$-element.
\end{lemma}
\begin{proof}
Letting $q$ be the order of a Sylow $p$-subgroup $S$ of $d(a)$,
we have $d(a)=d(a^q)\times S$ and thus $a H=s H$ for some $s\in S$.
\end{proof}

In the case that interests us, the group $H$  will contain no elements
of order $p$, and then the $p$-element in $a H\intersect d(a)$ will
be unique and of order $p$.

%End of file degenB02.tex
 %Overview
%degenB03.tex

\section{Minimization}

Lemma \ref{L:minimization} following will describe the situation
arising from consideration of a hypothetical minimal 
counterexample to Theorem \ref{thm:main} or \ref{thm:p-odd}.

\begin{definition}
Let $G$ be a group of finite Morley rank and $p$ a prime.
We say that $G$ is {\em $p$-degenerate} if $G$ contains
no infinite abelian $p$-subgroup.
\end{definition}

Observe that $G$ is $2$-degenerate if and only if $G$ is of degenerate
type. 

\begin{lemma}\label{L:minimization}
Let $G$ be a connected $p$-degenerate group of finite Morley rank
with nontrivial $p$-torsion, and of minimal Morley rank
among all such groups.
Then $\bar G=G/Z(G)$ is simple, and contains nontrivial $p$-torsion, while
no proper definable connected subgroup of $\bar G$ contains nontrivial
$p$-torsion. 
\end{lemma}
\begin{proof}
By our minimality hypothesis no proper connected subgroup of $G$
contains $p$-torsion. If $H<G$ is a nontrivial definable connected
normal subgroup, then passing to $G/H$ we contradict the minimality of
$G$. So $Z(G)$ is finite and $G/Z(G)$ is simple. It suffices
to show that $G/Z(G)$ contains nontrivial $p$-torsion.

Supposing the contrary, after passing to a quotient of $G$ we may
suppose that $Z(G)$ is a $p$-group.
We now introduce a function
$$\eta:G\to Z(G)$$
which though not a homomorphism will be covariant with respect to the
action of $Z(G)$.
This is defined as follows.

For $g\in G$, we consider the subgroup $\hat d(g)$, which
splits as $\hat d(g)^q\times S_g$, with $q$ the exponent of $Z(G)$ and
$S_g\le Z(G)$ the Sylow $p$-subgroup of $\hat d(g)$ (or of $d(g)$). 
So the projection $\pi_2:\hat d(g)\to S_g$ is well-defined, and we may
set $\eta(g)=\pi_2(g)\in Z(G)$. 

The desired covariance property is the following.
$$\hbox{$\eta(z g)=z\eta(g)$ for $z\in Z(G)$, $g\in G$}$$
Writing $g=g_0s$ with $g_0\in d(g)^q$ and $s\in S_g$,
we have $g_0^q\in d(z g)$ and as $d(g_0)$ is uniquely $p$-divisible
we have $g_0\in d(z g)$. But $z g=g_0zs$ and hence $zs\in d(z g)\le \hat
d(z g)$ as well, and our claim follows.

Now in view of the covariance of the map $\eta$, its fibers have
constant rank. Thus $G$ is partitioned by the fibers of $\eta$ into
finitely many sets of equal rank, and as $G$ is connected this yields
a contradiction.
\end{proof}

%End of file degenB03.tex
 %Definability
%degenB04.tex

\section{Genericity}

In the present section we suppose the following.

\begin{texteqn}{\dag}
$G$ is a $p$-degenerate simple group of finite Morley rank, $p$ is a prime,
$G$ contains nontrivial $p$-torsion, and no proper definable connected
  subgroup of $G$ contains nontrivial $p$-torsion.
\end{texteqn}

Let $q$ be a bound on the exponent of the $p$-torsion in $G$.

We will show that the generic elements of $G$ lie outside
every proper connected subgroup of $G$, and we will pin down their
location with sufficient precision to give useful structural
information. In particular we will show that the Sylow $p$-subgroup of
$G$ has exponent $p$, and that any two elements of order $p$ are
conjugate. 

We begin with a rank computation
similar to that of \cite[3.3]{ChJa-TMS}, but slightly more general.

\begin{lemma}\label{genericcovering}
Let $G$ be a group of finite Morley rank, $H$ a definable subgroup of $G$,
and $X\includedin G$ a definable set such that 
$$\rk(X\setminus \Union_{g\notin H} X^g)\ge \rk(H)$$
Then $\rk(\Union X^G)=\rk(G)$.
\end{lemma}
\begin{proof}
Let $Y=X\setminus \Union_{g\notin H} X^g$.
Let $\hat Y=\Union Y^H$. 
Then for $h\in H$ and $g\notin H$ we have 
$Y^h\intersect Y^g=(Y\intersect Y^{gh^{-1}})^h=\emptyset$, 
and thus $\hat Y\intersect \hat Y^g=\emptyset$. 
Furthermore $\rk(\hat Y)\ge \rk(Y)\ge \rk(H)$.

Now $H$ is the setwise stabilizer of $\hat Y$ in $G$
and thus
the family $\FF=\{\hat Y^g:g\in G\}$ has rank
$\rk(G)-\rk(H)$. As distinct elements of $\FF$ are
disjoint, we have
$$\rk(\Union \FF)=\rk(\hat Y)+\rk(\FF)\ge\rk(G)$$
Thus $\rk(\Union \hat Y^G)=\rk(G)$, and the same applies to 
$\Union X^G$. 
\end{proof}

In practice the group $H$ is the {\it generic stabilizer} of $X$ in
the following sense.

\begin{definition}
Let $G$ be a group of finite Morley rank and  $X\includedin G$
definable. 
The (generic) stabilizer $G[X]$ is the subgroup 
$$\{g\in G: \rk(X\triangle X^g)<\rk(X)\}$$
\end{definition}

The way we obtain the hypotheses of our lemma in
practice is by showing the following, in each instance.
$$\rk(X\setminus \Union_{g\notin G[X]}X^g)=\rk(X)\ge \rk(G[X])$$
We may then take $H=G[X]$.

Sometimes $X$ is even disjoint from its conjugates $X^g$ ($g\notin
G[X]$) but more often in applications this strong condition 
fails but our weaker condition holds.

\begin{lemma}\label{L:pgeneric}
Suppose $a$ is a $p$-element of $G$.
Then the set $\Union (a C\oo(a))^G$ is generic in $G$.
\end{lemma}

\begin{proof}
We remark that $C\oo(a)$ contains no $p$-torsion, by our hypothesis.

It will suffice to show that $a C\oo(a)$ is disjoint from its conjugates
$(a C\oo(a))^g$ for $g\notin C(a)$, as we may then apply Lemma
\ref{genericcovering}.  

So suppose we have an element $x\in (a C\oo(a))\intersect (a C\oo(a))^g$,
with $g\in G$.
Then $d(x)=A\oplus Z$ with $A$ $p$-divisible and $Z$ a cyclic
$p$-group, in fact $Z=\< a\>$. Accordingly $C\oo(Z)=C\oo(a)$
and $a$ is the only $p$-element in the coset $xC\oo(Z)$.
Similarly, since $x\in (a C\oo(a))^g$, 
we find that $a^g$ is the only $p$-element in $xC\oo(Z)$.

Thus $a=a^g$ and $g\in C(a)$.
\end{proof}

Now we can pass to structural consequences.
We emphasize that the clause $(\dag)$ has been assumed throughout.

\begin{lemma}\label{L:structure}
All nontrivial $p$-elements $a$ of $G$ 
are $G$-conjugate and of order $p$, and if $p=2$ then the
Sylow $2$-subgroup is elementary abelian.
\end{lemma}
\begin{proof}
For $a,b$ nontrivial $p$-elements of $G$, 
the sets $\Union (a C\oo(a))^G$ and $\Union
(b C\oo(b))^G$ are generic in $G$ and so we may suppose after conjugating
that $a C\oo(a)\intersect b C\oo(b)\ne \emptyset$. 
Now if $x\in aC\oo(a)\intersect bC\oo(b)$ then as in the proof of the 
previous lemma, $a=b$. 

So the $p$-elements are conjugate, and hence are all of order $p$. In
particular when $p=2$ the Sylow $2$-subgroup is elementary abelian.
\end{proof}

This is about as far as the genericity arguments take us, and it is
time to focus on the case $p=2$ and black box methods.

%End of file degenB04.tex
 %Minimization
%degenB05.tex

\section{Black box methods}

Now we specialize condition $(\dag)$ 
to the case $p=2$. Our hypotheses will therefore be as follows.

\begin{texteqn}[0.65]{\dag}
$G$ is a degenerate type simple group of finite Morley rank, 
$G$ contains involutions, and no proper definable connected
  subgroup of $G$ contains involutions
\end{texteqn}

For the moment we need not rely on the results of the previous
section. Rather we introduce a crucial case division, dispose of one
case using so-called black box group theoretic methods, and
then return to the other case using the information from the last
section. 

We turn to the case division.
Fix a conjugacy class of involutions $\CC$, and note that as this set
can be identified definably with $G/C(i)$ for any fixed $i\in \CC$, it
has Morley degree one. Thus the notion of ``generic element'' or pair
of elements in $\CC$ is robust. The first case we will treat is the
following one.

\begin{texteqn}[0.5]{\hbox{Case I}}
For generic and independent $i,j\in \CC$ the group $d(ij)$
contains no involution.
\end{texteqn}

The following two facts are elementary but important, and are used in
combination. 

\begin{lemma}\label{f:d(a)}
Let $G$ be a group of finite Morley rank and $a$ an element of $G$.
Suppose that the Sylow $2$-subgroups of $G$ have bounded exponent.
Then the following conditions on the group $d(a)$ are equivalent.
\begin{enumerate}
\item $d(a)$ contains no involutions.
\item $d(a)$ is $2$-divisible.
\item $d(a)$ is uniquely $2$-divisible.
\end{enumerate}
On the other hand, if $d(a)$ does contain an involution, then that
involution  is unique.
\end{lemma}

\begin{lemma}
Let $G$ be a group of finite Morley rank. 
Suppose that the Sylow $2$-subgroups of $G$ have bounded exponent.
Let $i,j$ involutions of $G$, and let $a=ij$. 
Then $d(i,j)=d(a)\sdprod \<i\>$, where $i$ acts by
inversion on $d(a)$. Furthermore, $i$ and $j$ are conjugate under the
action of $d(a)$ if and only if $d(a)$ contains no involution.
\end{lemma}

This is well known in the finite case and the infinite case is much
the same. The only noteworthy point is the fact that in the absence of
$2$-torsion the group is $2$-divisible, a point already used
repeatedly but relying essentially on the definability of the group in
question. 

This leads to consideration of the following partial functions from
the group $G$ to $C(i)$, for any fixed involution $i$, under the
hypothesis that the Sylow $2$-subgroups of $G$ have bounded exponent. 

\begin{enumerate}
\item $\zeta_0(g)$ is the unique involution in $d(i\cdot i^g)$, if
  $d(i\cdot i^g)$ contains an involution;
\item $\zeta_1(g)$ is the unique element in $gd(i\cdot i^g)\intersect
  C(i)$, if $d(i\cdot i^g)$ contains no involution.
\end{enumerate}

Indeed, under our hypothesis we have just seen that
$\zeta_0$ is well-defined on its domain. On the other hand,
if $d(i\cdot i^g)$ contains no involution then 
there is an element $x\in d(i\cdot i^g)$ conjugating $i$ to $i^g$,
so $gx^{-1}$ belongs to $C(i)\intersect gd(i\cdot i^g)$. 
As far as uniqueness is concerned, if $x,y\in C(i)\intersect gd(i\cdot
i^g)$, 
then $x^{-1}y\in C(i)\intersect d(i\cdot i^g)$ is both centralized and
inverted by $i$, hence an involution or trivial, and as we have ruled
out involutions we conclude $x=y$. One could also compute more directly
that $x^{-1}=x^i=x[x,i]=xi^gi$ and thus $x^2=ii^g$, so that $x$ is
uniquely determined within $d(i\cdot i^g)$, symbolically
$x=\sqrt{ii^g}$, with the square root operation restricted to
$d(i\cdot i^g)$, though this extra precision is useful mainly as a way
of verifying the existence of $x$.

The functions $\zeta_0$ and $\zeta_1$ are definable, 
because we can replace $d$ by $\hat d$
everywhere in their definitions, 
and the Sylow $2$-subgroups remain the same. The
uniqueness argument also relies on the properties of $\hat d$ given at
the outset, which mimic the properties of $d$.

\begin{lemma}\label{L:notI}
Case $(I)$ does not occur.
\end{lemma}
\begin{proof}
We fix $i\in \CC$ and consider the definable partial 
function $\zeta_1:G\to C(i)$ discussed above,
which is defined on a generic subset of $G$, namely
$$\zeta_1(g)\in C(i)\intersect gd(ii^g)$$
It follows by inspection of the definition that we have the covariance
property 
$$\zeta_1(cg)=c\zeta_1(g)$$
for $g$ in the domain of $\zeta_1$ and $c\in C(i)$. This implies that
the fibers of $\zeta_1$ are of constant rank, say $f$, and hence that
any subset of $C(i)$ of rank $r$ lifts under $\zeta_1$ to a subset of
$G$ of rank $r+f$. Now since $i\in C(i)\setminus C\oo(i)$, the group
$C(i)$ is disconnected and hence has disjoint subsets of full rank,
and these lift under $\zeta_1$ to disjoint generic subsets of $G$,
which contradicts the connectivity of $G$ \cite[Theorem 5.12]{BN}.
\end{proof}

%End of file degenB05.tex
 %Black box 
%degenB06.tex

\section{Proof of Theorem 1}

We suppose toward a contradiction that $G$ is a connected group of
degenerate type containing an involution.  Applying Lemma
\ref{L:minimization} we may suppose 
condition $(\dag)$ of \S5 holds:

\begin{texteqn}[0.65]{\dag}
$G$ is a degenerate type simple group of finite Morley rank, 
$G$ contains involutions, and no proper definable connected
  subgroup of $G$ contains involutions
\end{texteqn}

We fix a conjugacy class of involutions $\CC$ in $G$, and in view of 
Lemma \ref{L:notI} we suppose that the following holds.

\begin{texteqn}[0.5]{\hbox{Case II}}
For generic and independent $i,j\in \CC$ the group $d(ij)$
contains a unique involution.
\end{texteqn}

Applying Lemma \ref{L:structure} we find that the Sylow $2$-subgroups
of $G$ are elementary abelian. This then yields the following.

\begin{lemma}\label{Sec6:L1}
If $i,j$ are involutions and $k\in d(ij)$ is an involution,
then $i$ and $j$ are not conjugate under the action of $C(k)$.
\end{lemma}
\begin{proof}
We show first that $i$ and $ik$ are not conjugate in $C(k)$.
Suppose on the contrary $i^u=ik$ with $u\in C(k)$. 
Then $i^{u^2}=i$ and $u$ acts on the group $\<i,k\>$ 
as a nontrivial automorphism of order two. 
It follows that $d(u)$ contains a $2$-element with the same action, 
and as $G$ has abelian Sylow $2$-subgroups this is impossible.

On the other hand, as in the case of ordinary dihedral groups one may
see that the group $d(i,j)$ has two conjugacy classes of noncentral
involutions, represented by $i$ and $j$, and in particular $j$ is
conjugate to $ik$ under the action of $d(ij)$, and in particular
under $C(k)$. If $i$ is conjugate to $j$ under $C(k)$ then $i$ is
conjugate to $ik$ under $C(k)$ and we have a contradiction.
\end{proof}

Now we may conclude the proof of Theorem \ref{thm:main} 
by a model theoretic argument.

Fix a Sylow $2$-subgroup $S$ of $G$ and consider a pair
$i,j$ of involutions in $\CC$ which are independent and generic over
$S$, that is to say with the elements of $S$ treated as constants.
Define a subset $S_{i,j}\includedin S\times S$ as follows:
$$\{(s,t)\in S\times S: (i,k)\sim (s,t)\}$$
Here $k$ is the unique involution in $d(ij)$, and ``$\sim$'' refers
to conjugacy under the action of $G$. As $i$ and  $k$ commute,
the set $S_{i,j}$ is nonempty.

Now the pair $(i,j)$ and the pair $(j,i)$ have the same type over $S$,
so $S_{i,j}=S_{j,i}$.
As the involution $k$ is also the unique involution in $d(ji)$, this
means that $(i,k)$ and $(j,k)$ are conjugate to the same pairs in
$S\times S$, and hence to each other. But to conjugate $(i,k)$ to
$(j,k)$ in $G$ means that $i$ is conjugated to $j$ in $C(k)$. 
This contradicts the preceding lemma and completes the proof of
Theorem 1.

As mentioned in the introduction, we can also strengthen Theorem
\ref{thm:main} as follows.

\begin{namedtheorem}{Theorem \ref{thm:main2}}
Let $G$ be a connected group of finite Morley rank, and $i\in G$ 
an involution. 
Then the Sylow $2$-subgroup of $C(i)$ is infinite.
\end{namedtheorem}

\begin{proof}
Suppose toward a contradiction that $G$ and $i$ are a counter\-example
with $\rk(G)$ minimal. Then $i$ belongs to no proper definable
connected subgroup of $G$. 

By hypothesis $C\oo(i)$ is of degenerate type and hence contains
no involution. Hence for $c\in C\oo(i)$, $i$ is the unique involution
in $d(ci)$. 
Now for $g\in G$, if 
$(iC\oo(i))\intersect (iC\oo(i))^g\ne \emptyset$ 
and $x$ lies in the intersection, then $i,i^g\in d(x)$ and
thus $i=i^g$, and $g\in C(i)$. It follows that the distinct
conjugates $(iC\oo(i))^g$ are disjoint for $g\in G$,
and $iC\oo(i)$ has stabilizer $C(i)$, 

So by Lemma \ref{genericcovering} we have
$\rk(\Union[iC\oo(i)]^{G})=\rk(G)$.
In particular, the generic 
element of $G$ lies outside every proper definable connected subgroup.

We claim that a Sylow\oo\ $2$-subgroup $T$ of $G$ is a nontrivial
$2$-torus. It is nontrivial since $G$ is connected and contains an
involution, and it is a $2$-torus since otherwise it would have a
nontrivial unipotent factor $U$ which we may take to be normalized by
$i$, and then $C_U(i)$ would be infinite.  By \cite{Ch-GT} the generic
element of $G$ lies in a conjugate of $N\oo(T)=C\oo(T)$.
So if $C\oo(T)<G$ then the generic element lies
in a proper definable connected subgroup, a possibility we have just
ruled out. We are left with the possibility $C\oo(T)=G$, and then of
course $C_T(i)=T$.
So in any case we reach a contradiction.
\end{proof}

%End of file degenB06.tex
 %Main Theorem
%degenB07.tex

\section{Generic equations}

Now we take up Proposition \ref{prop:genericexponent}, concerning
generic equations of the form 
$$x^n=1\leqno(*)$$
In the  statement of that proposition it was assumed that $n$ was the
precise order of the generic element of $G$, but it will be convenient
here to work with the slightly broader condition $(*)$, so that the
order of a generic element is a fixed {\em divisor} of $n$.
As we have noted, Alt{\i}nel suggested that Theorem
\ref{thm:main} supplies the missing piece of the puzzle to treat
a substantial portion of this problem.

In the solvable case, with $n$ arbitrary, Jaber \cite{Ja-EG}
has shown
by an argument that uses results of Bryant and Wagner that when such
an equation holds generically then it holds everywhere \cite{Ja-EG}.

Let $G$ satisfy the equation $(*)$ generically and let $U$ be a
Sylow\oo\ 
$2$-subgroup of $G$. 
We break the proof into four steps. We will be arguing
inductively. 

\begin{enumerate}
\item $U$ is $2$-unipotent.
\item If $U$ is nontrivial then $G$ is not simple.
\item $U$ is normal in $G$.
\item $G=U\cdot C_G(U)$.
\end{enumerate}

For the first point, one may use the theory of
\cite{Ch-GT}, already cited above. If $G$ contains a nontrivial
decent torus in the sense of that reference, and if $T$ is a maximal
such, then $H=N\oo(T)$ 
is almost disjoint from its conjugates, that is
$H\setminus \Union_{g\notin N(H)}H^g$ is nongeneric in $H$,
and $\Union H^G$ is generic in $G$. It follows easily, under our
present hypotheses, that $H$ also satisfies the equation $(*)$
generically. 
But $H=C\oo(T)$ as well, and some coset of $T$ in $H$ must satisfy
$x^n=1$ generically, which is impossible. So there is no nontrivial
$p$-torus in $G$ for any $p$.
In particular, taking $p=2$, we find that $U$ is $2$-unipotent.

Now if $G$ is simple and $U$ is nontrivial, then $G$ is isomorphic to
an algebraic group, by the classification theorem for groups of even
type \cite{ABC:book}. But then $G$ contains a nontrivial
$p$-torus for almost all $p$, a contradiction. The same argument shows
that any definable simple section of $G$ must be of degenerate type.

Now we show that $U$ is normal in $G$. 

Using some general structure theory following from the classification 
of the simple groups of even type, one may argue as follows.
If $U$ is not normal in $G$, then the subgroup $B(G)$ generated
by all the unipotent 2-subgroups of $G$ is nonsolvable, and the
nonabelian factors of its socle are algebraic groups of even type; but
as we saw above, they must be of degenerate type. We may also argue
in a more {\it ad hoc} manner as follows.

We may suppose $U\ne 1$.
As $G$ is not simple, it is not definably simple, 
by a result of Zilber. It follows easily that $G$
contains a nontrivial definable connected proper normal subgroup, 
as otherwise consideration of $G/Z(G)$ (with $Z(G)$ finite) gives a
contradiction. 
Let $K$ be a minimal nontrivial definable connected normal
subgroup of $G$. Then $K$ is either abelian, or simple, and in the latter case 
is of degenerate type. 
In the abelian case, either $K\le U$ or $K$ is of degenerate type. 
The group $U$ acts on $K$, and if $K$ is of degenerate type this
action is trivial 
\cite[Proposition 2.8.6]{Alt-Hab},
\cite[Proposition 3.2]{AlCh-L*1}.

We may suppose inductively that $UK/K$ is normal in
$G/K$ and thus that $UK$ is normal in $G$. Either
$K\le U$, or $K$ is of degenerate type, in which case
$UK$ splits as $U\times K$ since $U$ centralizes $K$.
In either case it follows at once that $U$ is normal in $G$.

Finally, take a maximal connected definable $G$-invariant
series 
$$1=U_0<U_1<\dots <U_n=U$$ 
for $U$ under the action of $G$.
As $U$ acts trivially on each successive factor $V_i=U_i/U_{i-1}$, 
$\bar G=G/U$ acts on each factor as a
group of degenerate type. 
By \cite[Lemma 3.13]{AlCh-LSE}, the Borel subgroups
of $\bar G/C_{\bar G}{V_i}$ are good tori, hence trivial by our current
assumptions. It follows that $G$ acts trivially on each such factor,
and thus every definable
$2^\perp$-subgroup of $G$ (including the finite ones) acts trivially
on $U$. But $G$ is generated by $U$ and its $2^\perp$-subgroups, as it
suffices to look at $d(g)$ as $g$ runs over $G$, and thus 
$G=U\cdot C_G(U)$, as claimed.

With this analysis in hand, 
we may consider the quotient $G/U$,
which again satisfies $(*)$ generically. As $G/U$ is of degenerate
type it contains no involutions, and thus we may replace $n$ by its
odd part. Note that the precise order of the generic element may
change in passing from $G$ to $G/U$.

%End of file degenB07.tex
 %odd primes
%degenB08.tex

\section{The case of odd primes}

We now prove Theorem \ref{thm:p-odd}. This reads as follows.

\begin{namedtheorem}{Theorem \protect\ref{thm:p-odd}}
Let $G$ be a connected group of finite Morley rank
whose Sylow $p$-subgroups are finite. Then $G$ contains no elements of
order $p$.
\end{namedtheorem}

Here, as it happens, the two possible notions of Sylow $p$-group
(maximal $p$-group or maximal locally nilpotent $p$-group)
give the statement in question two possible meanings, and we prove
this actually in the stronger of the two possible forms, 
which may be put as follows.
$$\hbox{A connected $p$-degenerate group of finite Morley rank has no
$p$-torsion.}$$ 

\begin{proofof}{Theorem \ref{thm:p-odd}}
Suppose toward a contradiction that $G$ is a $p$-degenerate group 
of finite Morley rank with nontrivial $p$-torsion.
Applying Lemma \ref{L:minimization} we
may suppose that $G$ is simple and that no proper definable connected
subgroup of $G$ contains nontrivial $p$-torsion.

By Theorem \ref{thm:main}, the prime $p$ is odd.

Now applying Lemma \ref{L:structure}, it follows that all elements of order
$p$ are conjugate in $G$, and in particular each such element is
conjugate to its inverse. Hence $G$ contains involutions, 
and is not of degenerate type.

If $G$ contains a nontrivial unipotent $2$-subgroup, then by the
results of \cite{ABC:book}, $G$ is a simple algebraic group over an
algebraically closed field of characteristic two. 
But then $G$ contains unbounded $p$-torsion for all odd primes $p$.

It follows that the Sylow\oo\ $2$-subgroup of $G$ is a nontrivial
$2$-torus. In particular $G$ contains nontrivial ``decent tori'' in
the sense  of \cite{Ch-GT}, and by the analysis given at the end of
that reference, if $T$ is a maximal decent torus (that is, maximal
among definable abelian subgroups of the form $d(T_0)$ with
$T_0$ a divisible torsion subgroup), then
$\Union (N\oo(T))^G$ is generic in $G$.

But by Lemma \ref{L:pgeneric}, the generic element of $G$ lies outside
every proper definable connected subgroup of $G$, and we have a
contradiction.
\end{proofof}

As an application, we can prove the Genericity Conjecture
\ref{con:genericity} for minimal
connected simple groups.

\begin{proposition}
Let $G$ be a minimal connected simple group of finite Morley rank.
Then the set of elements of $G$ which belong to some connected
nilpotent subgroup of $G$ contains a definable generic subset of $G$.
\end{proposition}
\begin{proof}

Assume the contrary. 
Then the group contains no decent torus, in view
of the results of \cite{Ch-GT} and a result of Fr\'econ
\cite[3.5]{ChJa-TMS}.
So $G$ is of degenerate type (making use of \cite{ABC:book} to eliminate
even and mixed type).
On the other hand, there must be $p$-torsion for some prime $p$, so
$G$ is not $p$-degenerate. So there is an infinite abelian
$p$-subgroup of $G$, and as there is no decent torus there 
is a nontrivial $p$-unipotent subgroup in $G$. 

Let $B$ be a Borel subgroup of $G$ with $U_p(B)$ nontrivial.
As $B$ contains no decent torus, $B/U_p(B)$ contains no $p$-torsion,
and $B=U_p(B)C_B(U_p(B))$.

A generic element $g$ of $B$ has the property that $d(g)$ meets
$U_p(B)$; otherwise, we would have a generic subset
$X$ of $B$ such that for $g\in X$ the group $d(g)$ is $p$-torsion
free, and then multiplying $X$ by a nontrivial element of $Z(U_p(B))$ yields
a contradiction.

It follows that $B$ is generically disjoint from its conjugates 
as otherwise we would have distinct Borel subgroups meeting
$U_p(B)$ nontrivially, contradicting \cite[2.1]{Bu-BM}.

So a generic element of $G$ is conjugate to a generic element of $B$,
which lies in a Carter subgroup of $B$ by a result of Fr\'econ
\cite[3.5]{ChJa-TMS}.
\end{proof}

\begin{namedtheorem}{Theorem \protect\ref{thm:p-odd2}}
Let $G$ be a connected group of finite Morley rank and $a\in G$ a
$p$-element.
Then $C\oo(a)$ contains an infinite abelian $p$-subgroup.
\end{namedtheorem}
\begin{proof}
We follow closely the proof of Theorem \ref{thm:main2}, analyzing a
minimal counterexample.

First one argues that $\Union [aC\oo(a)]^G$ is generic in $G$.  This
proceeds as before. One concludes that the generic element of $G$ lies
outside every definable proper connected subgroup.  
 
We claim that $G$ contains no connected normal abelian
definable subgroup $A$. Supposing the contrary, by induction
there is an infinite connected abelian $p$-group $\bar P$ contained in
$C_{G/A}(a)$, with preimage $P$ in $G$. 

In particular $d(P)\< a\>$ is solvable. By the Hall theory for solvable
groups of finite Morley rank $a$ normalizes a Sylow\oo\ $p$-subgroup $Q$ of
$d(P)$, and $Q$ covers $\bar P$.
But then $[a,Q]\le Q\intersect A$ and as $Q>Q\intersect A$ we find
$C_Q(a)$ is infinite and reach a contradiction.

So $G$ contains no connected normal abelian definable subgroup and in
particular $Z(G)$ is finite.
It now follows that there are no nontrivial $\ell$-tori
in $G$ for any $\ell$: if $T$ is a nontrivial $\ell$-torus
then the generic element lies in a conjugate of $C\oo(T)$,
and as $T$ cannot be central in $G$ this is a contradiction.

Now $\Union [a^{-1}C\oo(a)]^G$ is also generic in $G$, and it follows
easily 
as in the proof of Lemma \ref{L:structure}
that $a$ and $a^{-1}$ are conjugate in $G$. This implies that $G$
contains an involution. As $G$ contains no nontrivial $2$-torus, it is
of even type.

But then the group $B(G)$ generated by the unipotent
$2$-subgroups of $G$ is a $K$-group by \cite[Proposition 3.4]{AlCh-L*1} and
the classification of simple groups of even type, and thus has no
definable connected simple sections of degenerate type.
As there are no definable simple sections of even type either,
the group $B(G)$ must be a solvable connected normal subgroup of $G$, and $G$
contains an abelian connected normal subgroup, a contradiction.
\end{proof}

%End of file degenB08.tex
 %Generation
%degenB09.tex

\section{Generation}

We turn to Theorem \ref{thm:generation}, or rather the following,
which is marginally stronger in view of
Theorem \ref{thm:main} and Lemma \ref{lem:C(i)connected} below.

\begin{proposition}\label{prop:generation2}
Let $H$ be a group of finite Morley rank without involutions,
and $V$ a $4$-group acting
definably on $H$. Then 
$$H=\<C(v):v\in V^\#\>$$
\end{proposition}

\begin{proof}
Work in $G=H\cdot V$, with $V=\<i,j\>$. 

Fix $h\in H$.
Then $i$ and  $j^h$ are not conjugate in $G$, so there is some
$u\in I(G)$ commuting with both (indeed, $u\in d(i,j^h)$).

As $V$ is a Sylow 2-subgroup of $C(i)$ and $u\in C(i)$, 
there is some $h_0\in C(i)$ so that
$u^{h_0}\in V$. Similarly there is $h_1\in C(u^{h_0})$ so that
$(j^{hh_0})^{h_1}\in V$. Then $j^{hh_0h_1}=j$ and $hh_0h_1\in C(j)$, 
so $h\in C(j)C(u^{h_0})C(i)$. This proves the claim.
\end{proof}

Now in order to recover the statement of Theorem \ref{thm:generation},
we apply Theorem \ref{thm:main} and also the following general lemma.

\begin{lemma}\label{lem:C(i)connected}
Let $G$ be a connected group of finite Morley rank without
involutions, and $i$ an involutory automorphism of $G$. Then $C_G(i)$
is connected.
\end{lemma}
\begin{proof}
We use a ``black box'' argument in the style of \S5.

We consider the group $\hat G=G\sdprod \<i\>$. For any pair of
involutions $i,j\in \hat G$ we have $ij\in G$ and thus $d(ij)$
is $2$-torsion-free. It follows that the map
$$\zeta_1:\hat G\to C_{\hat G}(i)$$
introduced in \S5 is well-defined.
Observe that the restriction of $\zeta_1$ to $G$ carries
$G$ into $C_G(i)$, and is $C_G(i)$-covariant. Thus we find as usual 
$\deg(C_G(i))\le \deg(G)=1$ and  $C_G(i)$ is connected.
\end{proof}

We also give a ``lifting'' lemma, in the spirit of extending as much
as possible of the solvable theory to the general case.

\begin{lemma}\label{L:C-lift}
Let $G$ be a group of finite Morley rank
and $i$ an involution acting on $G$. 
Let $H\normalin G$ be definable and $i$-invariant, without
involutions, and with $i$ acting
trivially on $G/H$. Then $G=H\cdot C(i)$.
\end{lemma}
\begin{proof}
Let $g\in G$, and set $h=[i,g]\in H$. 
Then $i$ inverts $h$.

As $H$ contains no involutions, the group $d(h)$ is $2$-divisible.
Take $h_1\in d(h)$ with $h_1^2=h$. Then $i$ inverts $h_1$ and hence
$$[i,h_1]=h_1^2=[i,g]$$
Thus $gh_1^{-1}\in C(i)$, and $g\in C(i)H=HC(i)$.
\end{proof}

%End of file degenB09.tex
 %Generic equations
%degenB10.tex

\section{Afterword: Black-box Groups}

The methods used for the proof of the main result are relatively
self-contained. One of the main ingredients comes from an unusual source:
``black box group theory'' \cite{KS-BB}, via a line 
of thought represented by \cite{Br-CI,PW-BB,altseimer-borovik,Bo-CI}.
The subject as a whole deals with the problem of
computing in large finite groups which are given in such a form that
one can extract elements {\em at random,} 
and perform limited operations or tests on these
elements. Among the problems in this area are the determination as to
whether the group in question is simple, and its identification if it
is. The issue of identification of black box groups
is a subject which has remarkable affinities with the subject
of groups of finite Morley rank, which can be traced back to a
preoccupation with ``generic'' elements. In the probabilistic setting,
this refers to the kinds of elements that appear with probability 1,
while in the model theoretic setting this refers to the sets of
elements which have maximal dimension (Morley rank). Of course, here
one may only consider sets which are either measurable or definable,
respectively. In black box group theory one can fall back on the
classification of the finite simple groups \cite{GLS-SG}, 
whereas in groups of
finite Morley rank this is the problem which is under investigation.
However the analogy can be maintained, because with or without
a classification theorem, the problem is one of {\em recognition}
of the specific group with which one is presented, by methods 
allowed by the corresponding framework. The difference is in outcome:
black box group theory delivers practical algorithms (implemented, for
example, in GAP) while the theory of groups of finite Morley rank is a
conventional mathematical theory dealing in theorems.

In either case, at a certain point, just as in the case of
conventional group theory, one requires information about 
{\em centralizers of involutions}. Heretofore in dealing with groups of
finite Morley rank we have followed the lead of the finite group
theorists, who have a powerful and elegant range of techniques for
dealing with this problem. But there is another very direct way to
gain a measure of control over the centralizer of an involution which
appears to have surfaced first within black box group theory.
It is based on elementary properties of dihedral groups which
are essential to classical fusion analysis, but via a certain (partial)
function $\zeta$ from the group $G$ to the centralizer of an
involution which appears to be entirely useless in conventional group
theory, probably because it is a partial function, and comes
into its own only when it is {\em generically defined} (that is,
its domain is a {\em generic subset} of $G$). This is the function
$\zeta_1$ which we encountered in \S5, along with its companion
$\zeta_0$.

From a technical point
of view, the virtue of the function $\zeta$ is that it preserves
{\em uniform distribution} in the probabilistic setting and 
{\em connectivity} in the finite Morley rank setting. Accordingly,
if the function $\zeta$ is generically defined, then in the black box
setting one can deduce that the centralizer of an involution is again
a black box group, and in the finite Morley rank setting one can
deduce that the centralizer of an involution is connected if the
original group was. We insist here on these parallels because they
appear to go more than skin deep. 

We mention for the sake of finite group theorists that this technique
produces a version of the celebrated $Z^*$-theorem \cite{Gl-Z*}, 
proved by methods
that have no known analog in the finite case. Indeed, our version
assumes connectivity of the ambient group.

\begin{theorem}[$Z^*$]\label{Z*}
Let $G$ be a connected group of finite Morley rank, $S$ a Sylow
$2$-subgroup in $G$ and $i \in S$ an involution.
Then either

\begin{itemize}
\item[(a)] $i$ is conjugate in $G$ to another involution in $S$, or

\item[(b)] $C_G(i)$ is connected.
\end{itemize}
\end{theorem}
\begin{proof}

Let $\CC$ be the conjugacy class $i^G$. 
Suppose condition $(a)$ fails: $\CC\intersect S=\{i\}$. 
We may suppose $i\notin Z(G)$.

If $j\in \CC$ commutes with $i$, then
$j=i$. 
Otherwise, the group $\<i,j\>$ could be conjugated into $S$,
giving two conjugates of $i$ in $S$, at least one of which is not
$i$.

Fix an involution $j\in \CC$.

If $d(ij)$ contains an element $t$ of order $4$, then $i$ inverts $t$
and the subgroup $\<i,t\>$ can be conjugated into $S$. 
As $\CC\intersect S=\{i\}$ it follows that $i$ inverts an element
$t'\in S$ of order $4$, and then $i^{t'}\ne i$, contradicting our
assumptions.  

So $d(ij)$ contains no element of order $4$.
Now suppose $d(ij)$ contains an involution $k$. Then (much as in the
finite case) $j$ is conjugate to $ik$, and hence $i$ is conjugate to
$ik$, which commutes with $i$, giving a contradiction.

So $d(ij)$ contains no involutions. 
Now holding $i$ fixed,
we need to vary $j$ and to consider $d(ij)$ 
as a function of the element $j$, and
we need a generically definable function here.
Let $j$ be generic in $\CC$ and independent from $i$,
and let $\phi(x,y)$ be a formula such that
$\phi(x,ij)$ defines $d(ij)$. The set of involutions $j'\in \CC$ such that
$\phi(x,ij')$ defines an abelian group containing $ij'$, inverted by
$i$, and without involutions, is a generic subset of $\CC$.
So letting $\hat d(ij')$ be the group defined by $\phi(x,ij')$ for
such $j'$, we can use $\hat d$ as a definable approximation to $d$
and define a covariant function $\zeta_1:G\to C(i)$ as we did earlier, 
then deduce from the connectedness of $G$ that  $C(i)$ is connected.
\end{proof}

The methods coming from black box group theory offer the 
outstanding advantage that they do not rely much on induction.
In this they resemble the transfer method, which however does not have
a direct analog in our context. As degenerate type groups are very poorly
understood, this is an essential point. We do in fact use an inductive
argument, but only with respect to the presence of involutions; we 
do not assume anything else about the simple sections which may be
present in the group.

The arguments given here emerged gradually, and have been considerably
simplified over time. Earlier arguments used the Alperin-Goldschmidt
Theorem \cite{Co-AG} and some of the $0$-unipotence theory developed
by the second author \cite{Bu-SF,Bu-S0}. As these give only special
cases of the results given here, we will not elaborate on this point,
but as there is a good deal more to be done in the study of groups of
degenerate type, these alternative techniques may yet have a role to play.

%End of file degenB10.tex
 %Afterword: generation, applications, etc.

\bibliographystyle{plain}
%\bibliographystyle{alpha}
%invol.bbl

% End of file invol.bbl

\medskip
\end{document}